\documentclass[12pt,twoside]{article}
\usepackage{amssymb,hyperref}
\usepackage{latexsym}
\usepackage{amsmath}
\usepackage[all]{xy}
\usepackage{sectsty}
\usepackage[left=2cm,top=2.5cm,bottom=3cm,right=3cm]{geometry}
\usepackage{fancyhdr}
\sectionfont{\fontfamily{ptm}\selectfont}
\allsectionsfont{\mdseries \centering \scshape}
\makeatletter
\def\@seccntformat#1{\csname the#1\endcsname.\quad}
\makeatother
\newtheorem{theorem}{Theorem}[section]
\newtheorem{proposition}[theorem]{Proposition}
\newtheorem{lemma}[theorem]{Lemma}

\newtheorem{definition}[theorem]{Definition}
\newtheorem{remark}[theorem]{Remark}
\newenvironment{proof}{\noindent {\bf Proof}~} {\hfill $\Box$\medskip}

\newcommand{\ox}{{\cal O}_{X} }
\newcommand{\oc}{{\cal O}_{C} }
\newcommand{\op}{{\cal O}_{\mathbb{P}^N_k}(1) }
\newcommand{\oxh}{{\cal O}_{X_h} }
\begin{document}
\title{Some Koszul Rings from Geometry}
\author{Krishna Hanumanthu}
\date{}
\maketitle
\begin{abstract}
We give examples of Koszul rings that arise naturally in algebraic geometry. In the first part, we prove a general result on Koszul property associated to an ample line bundle on a projective variety. Specifically, we show how Koszul property of multiples of a base point free ample line bundle depends on its Castelnuovo-Mumford regularity. In the second part, we give examples of Koszul rings that come from adjoint line bundles on minimal irregular surfaces of general type. 
\end{abstract}

\section*{Introduction}
Let $k$ be  a field. A standard graded $k$-algebra $R = k \bigoplus R_1 \bigoplus R_2 ...$ is  said to be {\it Koszul} if $k$ has a linear minimal resolution as an $R$-algebra.

Let $$...\rightarrow E_p \rightarrow E_{p-1} \rightarrow ... \rightarrow E_1 \rightarrow E_0 \rightarrow k \rightarrow 0$$ be a minimal resolution of $k$ over $R$. Then R is Koszul if and only if $E_0 = R$ and $E_p = {R(-p)}^{\oplus r(p)}$ for $p \geq 1$. Equivalently, Tor$_i^R(k,k)$ has pure degree $i$ for all $i$.

Koszul algebras were introduced by Stewart Priddy \cite{10} and they have applications in many areas of mathematics, such as algebraic geometry, commutative algebra and representation theory to name a few. For a sample of these applications, see \cite{11}, \cite{12}, \cite{13}, \cite{14}. See \cite{15} for a general introduction to Koszul property with historical notes. \cite{14} also has a general treatment of Koszul property. 

Part of the algebraic geometer's interest in Koszul rings stems from the following observation: 

Let $L$ be a very ample line bundle on a projective variety $X$ over $k$. Let $I_X$ be the ideal defining $X$ under the embedding in a projective space defined by $L$. Define 
$$R(L) = \bigoplus_{n=0}^{\infty} H^0(X, L^{\otimes n}).$$  

If $R(L)$ is Koszul then $X$ is projectively normal and $I_X$ is generated by quadrics. In the notation of $N_p$ property, this means that $L$ satisfies the property $N_1$. (See \cite{8}, 1.8.D for details on $N_p$ property.) If $R(L)$ is Koszul we say that $L$ has Koszul property.

There are several results establishing Koszul property for line bundles on curves. For instance,  see \cite{17}, \cite{18}, \cite{20} and \cite{19}. Koszul property for line bundles on elliptic ruled surfaces is studied in \cite{1}. Koszul property for adjoint line bundles on regular surfaces is studied in \cite{2}. Some general results on Koszul property for adjunction bundles is discussed in \cite{3}. 

Quite generally, high enough powers of ample line bundles have Koszul property (\cite{16}). The relation between the precise powers that achieve Koszul property  and Castelnuovo-Mumford regularity of the bundle is of general interest. In the first part of this paper (Section \ref{part1}), we prove a general result (Theorem \ref{final}) establishing such a relation. If $B$ is a base point free  ample 
line bundle on a projective variety $X$ and if $reg(B)$ is $r$ (cf. Definition \ref{new}), we show that $B^{\otimes n}$ has Koszul property for $n \geq r+1$.  This result is also proved  in a preprint (\cite{22}) invoking the notion of multigraded regularity. The proof in \cite{22} and our proof are both motivated by Theorem 1.3 in \cite{2} and essentially follow the methods developed there. 

A similar result is proved in \cite{16} and \cite{23}\footnote{I sincerely thank Burt Totaro for bringing this result to my notice.}. Let $R$ be a polynomial ring and let $I \subset R$ be a homogeneous ideal. Set $A = R/I$. In these papers, authors develop useful criteria to determine if the $d$-th Veronese subring $A_{(d)}$ of $A$ is Koszul. One of their results says that if $d \geq reg(I)/2$, then $A_{(d)}$ is Koszul. In our situation, this means the following: let $B$ be a base point free, ample line bundle that defines a map whose image is projectively normal in the projective space. Then $B^{\otimes d}$ has Koszul property for $d \geq reg(B)/2$. 

In the second part of the paper (Section \ref{part2}) we give examples of Koszul rings associated to certain adjoint line bundles on a minimal irregular surface of general type. This extends an analogous result for regular surfaces in \cite{2}.  As mentioned above, the Koszul property implies $N_1$, but the converse is, in general, not true (\cite{21}). Our theorem establishes the converse in this case. Our method is similar to \cite{4}.

Establishing $N_1$ property involves proving that a certain multiplication map of global sections of vector bundles is surjective. Koszul property is equivalent to the surjectivity of infinitely many multiplication maps of global sections of certain vector bundles, first of which is the multiplication map that appears in the $N_1$ property. In most examples of Koszul rings arising in algebraic geometry, the surjectivity required for Koszul property is proved by methods very similar to those used in establishing the $N_1$ property, after an appropriate inductive framework is set up. However, in the case of adjoint line bundles on irregular surfaces that we study, the methods used in establishing the $N_1$ property (\cite{4}) do not work for the subsequent surjections required for Koszul property. This suggests a potential example where $N_1$ property does not imply Koszul property. In this paper, we establish the Koszul property under a stronger assumption than was made in \cite{4}, namely the canonical bundle is base point free.


Many of the results cited here are directly connected to the cases we study. They represent only a sliver of the research on Koszul property in algebraic geometry. There are many results of a similar flavor that we do not mention here, but that are interesting nevertheless. 

{\bf Acknowledgement}:  B.P. Purnaprajna introduced me to this subject, taught me the key concepts and guided me throughout this work. I thank him for this and for his continued encouragement.



\section{Preliminaries}\label{section1}
Let $k$ be a field and let $X$ be a projective variety over $k$.  

{\bf Notation}: For a coherent sheaf $F$ on $X$, we write $H^i(F)$ to denote the $i$th sheaf cohomology group 
$H^i(X, F)$. 

Let $L$ be a line bundle on $X$. As before consider the ring:
$$R(L) = \bigoplus_{n=0}^{\infty} H^0(L^{\otimes n}).$$  

The question of whether $R(L)$ is Koszul has a nice cohomological interpretation due to Lazarsfeld. 


Given any vector bundle $F$ on $X$ that is generated by its global sections, we have a canonical surjective map\\
\begin{eqnarray}
H^0(F) \otimes {\cal O}_{X} \rightarrow F
\end{eqnarray}

Let $M_F$ be the kernel of this map. We have then the natural exact sequence
\begin{eqnarray}
0 \rightarrow M_F \rightarrow H^0(F) \otimes {\cal O}_{X} \rightarrow F \rightarrow 0
\end{eqnarray}

Now set $M^{(0),L} := L$.

If $L$ is globally generated, define 
$$M^{(1),L}: = M_L \otimes L = M_{M^{(0),L}} \otimes L.$$

If $M^{(1),L}$ is generated by its global sections, define 
 
$$M^{(2),L}: = M_{M^{(1),L}} \otimes L.$$

Inductively, define 
$M^{(h),L}: = M_{M^{(h-1),L}} \otimes L$, provided that $ M_{M^{(h-1),L}}$ is generated by its global sections. 

Then we have the following proposition that characterizes the Koszul property of $L$ in terms of certain cohomology groups. 

\begin{proposition}\label{p}\textup{[\cite{3}, Lemma 1]}
Let $X$ be a projective variety over a
field $k$. Assume that $L$ is a base-point-free line bundle on $X$ such that the 
vector bundles $M^{(h),L}$ are globally generated for every $h \geq 0$. If 
$H^1(M^{(h),L} \otimes L^s) = 0$
for every $h \geq 0$ and every $s \geq 0$, then $R(L)$ is Koszul.
\end{proposition}

\section{Preparatory lemmas}
In this section we will list and prove some well known results that will be used repeatedly in what follows. $k$ is any field and $X$ is a projective variety over $k$.

\begin{lemma}\label{horace}
Let $E$ and $L_1, L_2,...,L_r$ be coherent sheaves on $X$. Consider the multiplication maps

$\psi: H^0(E) \otimes H^0(L_1 \otimes ...\otimes L_r) \rightarrow H^0(E \otimes L_1 \otimes ... \otimes L_r)$,

$\alpha_1: H^0(E) \otimes H^0(L_1) \rightarrow H^0(E \otimes L_1)$,

$\alpha_2: H^0(E \otimes L_1) \otimes H^0(L_2) \rightarrow H^0(E \otimes L_1 \otimes L_2)$,

...,

$\alpha_r: H^0(E\otimes L_1 \otimes ... \otimes L_{r-1}) \otimes H^0(L_r) \rightarrow H^0(E \otimes L_1 \otimes ... \otimes L_{r})$.

If $\alpha_1$,...,$\alpha_{r-1}$ are surjective then so is $\psi$.
\end{lemma}
\begin{proof}
We have the following commutative diagram where $id$ denotes the identity morphism:
\begin{displaymath}
\xymatrix @R=2pc @C=3pc{
H^0(E) \otimes H^0(L_1) \otimes ... \otimes H^0(L_r) \ar[r]^{\alpha_1 \otimes id} \ar[d]^{\phi} & H^0(E \otimes L_1) \otimes H^0(L_2) \otimes ... \otimes H^0(L_r) \ar[d]^{\alpha_2 \otimes id} \\
H^0(E) \otimes H^0(L_1 \otimes ... \otimes L_r)    \ar[dd]^{\psi} & H^0(E \otimes L_1 \otimes L_2) \otimes H^0(L_3) \otimes ... \otimes H^0(L_r)  \ar[d]^{\alpha_3 \otimes id} \\
   &...\ar[d]^{\alpha_{r-1} \otimes id}\\
H^0(E \otimes L_1 \otimes ... \otimes L_r)   &\ar[l]_{\alpha_r} H^0(E \otimes L_1 \otimes ... \otimes L_{r-1}) \otimes H^0(L_r)}
\end{displaymath}

Since $\alpha_1, \alpha_2,...,\alpha_r$ are surjective and this diagram is commutative, a simple diagram chase shows that  $\psi$ is surjective.
\end{proof}

\begin{lemma}\label{globalgen}
Let $F$ be a locally free sheaf and $A$ an ample line bundle on $X$. If the multiplication map
$H^0(F \otimes A^{\otimes n}) \otimes H^0(A) \rightarrow H^0(F \otimes A^{\otimes n+1})$ is surjective for every $n \geq 0$, then $F$ is generated by its global sections. 
\end{lemma}
\begin{proof}
Since $A$ is ample, there exists  $m \geq 0$ such that $F \otimes A^{\otimes m}$ is generated by global sections. In other words, the morphism of sheaves
$\nu: H^0(F \otimes A^{\otimes m}) \otimes \ox \rightarrow F \otimes A^{\otimes m}$ is surjective. 

The hypothesis implies, by Lemma \ref{horace}, that $\psi: H^0(F) \otimes H^0(A^{\otimes m}) \rightarrow H^0(F \otimes A^{\otimes m})$ is surjective. 

Consider now the commutative diagram:
\begin{displaymath}
\xymatrix @R=2pc @C=3pc{
H^0(F \otimes A^{\otimes m}) \otimes \ox \ar[rd]^{\nu} &  \\
H^0(F) \otimes H^0(A^{\otimes m}) \otimes \ox \ar[u]^{\psi \otimes id}  \ar[d]_{id \otimes \phi} & F \otimes A^{\otimes m}\\
H^0(F) \otimes A^{\otimes m} \ar[ru]^{\mu}
}
\end{displaymath}

Since $\psi \otimes id$ and $\nu$ are surjective, a diagram chase shows that $\mu: H^0(F) \otimes A^{\otimes m} \rightarrow F \otimes A^{\otimes m}$ is surjective. 

As $A^{\otimes m}$ is an invertible sheaf, the surjectivity of $\mu$ shows that $F$ is generated by global sections. 
\end{proof}

\begin{lemma}\label{key} \textup{[CM Lemma, \cite{5}]}
Let $E$ be a base-point free line bundle on $X$ and let $F$ be a coherent sheaf on $X$. If $H^i(F \otimes E^{-i}) = 0$ for $i \geq 1$, then the multiplication map 
$$H^i(F \otimes E^{\otimes i}) \otimes H^0(E) \rightarrow H^i(F \otimes E^{\otimes i+1})$$
is surjective for all $i \geq 0$. 
\end{lemma}

Let $N$ be a globally generated vector bundle and let $A$ be a line bundle on $X$.

We have a short exact sequence
\begin{eqnarray}
0 \rightarrow M_N \rightarrow H^0(N) \otimes \ox \rightarrow N \rightarrow 0 \label{n}
\end{eqnarray}

\begin{remark}\label{-1}
$H^1(M_N \otimes A) = 0$ if the following two conditions hold.
\begin{itemize}
\item The multiplication map $H^0(N) \otimes H^0(A) \rightarrow H^0(N \otimes A)$ is surjective.
\item $H^1(A) = 0$.
\end{itemize}

This is easy to see: tensor the sequence (\ref{n}) by $A$ and take global sections:

$..\rightarrow H^0(N) \otimes H^0(A) \rightarrow H^0(N \otimes A) \rightarrow H^1(M_N \otimes A)
\rightarrow H^0(N) \otimes H^1(A) \rightarrow ...$
\end{remark}
\begin{remark}\label{-2}
$H^2(M_N \otimes A) = 0$ if the following two conditions hold.
\begin{itemize}
\item $H^1(N \otimes A) = 0$.
\item $H^2(A) = 0$.
\end{itemize}
This is easy to see: tensor the sequence (\ref{n}) by $A$ and take global sections:

$..\rightarrow H^1(N \otimes A) \rightarrow H^2(M_N \otimes A)
\rightarrow H^0(N) \otimes H^2(A) \rightarrow ...$
\end{remark}

\section{Koszul ring associated to an ample line bundle on a projective variety}\label{part1}
In this section $X$ denotes an arbitrary projective variety over a field $k$.

Let $B$ be a base point free ample line bundle on $X$. 

\begin{definition} \label{new} \textup{[\cite{8}, Definition 1.8.4.]} Let $m \geq 0$. We say that $B$ is $m-$regular (with respect to $B$) if 
$$H^i(B^{\otimes m+1-i}) = 0, \text{~for~} i>0.$$
\end{definition}


If $B$ is $m$-regular, then it is $(m+1)$-regular (\cite{8}, Theorem 1.8.5.(iii)). We define the {\it regularity} of $B$ to be $m$ if $B$ is $m$-regular, but {\it not} $(m-1)$-regular.  

This notion of regularity is related to the classical notion of Castelnuovo-Mumford regularity as follows:

Let $f: X \rightarrow \mathbb{P}^N_k$ be the morphism to a projective space defined by $B$. Note that such a morphism exists because $B$ is base point free. Let $L = f_{\star}(B)$. 

In the classical setting, we say $L$ is $m$-regular if $H^i\big{(}\mathbb{P}^N_k, L(m-i)\big{)} = 0$ for $i \geq 0$.

Since $f^{\star}(\op) = B$, by the projection formula we get
$$f_{\star}(B^{\otimes 2}) = f_{\star}(B \otimes f^{\star}(\op)) 
= f_{\star}(B) \otimes \op = L(1).$$

By induction, we obtain for any $r \geq 1$,
\begin{eqnarray}\label{reg}
f_{\star}(B^{\otimes r}) = L(r-1)
\end{eqnarray}

Since the morphism $f$ is finite, we have $H^i(X, A) \cong H^i(\mathbb{P}^N_k, f_{\star}(A))$ for any 
sheaf $A$ on $X$. 

Hence, by (\ref{reg}), $B$ is $m$-regular in the sense of Definition \ref{new} if and only if $f_{\star}(B)$ is $m$-regular in the sense of Castelnuovo-Mumford.

Suppose now that $B$ is $(r-1)$-regular (in the sense of Definition \ref{new}). Then since $B$ is $(n-1)$-regular for all $n \geq r$, we have   
\begin{eqnarray}\label{cm}
H^i(B^{\otimes n-i}) = 0 ~\text {for all}~ i \geq 1 ~\text {and}~  n \geq r
\end{eqnarray}

Set $L = B^{\otimes r}$. We prove that $R(L)$ is a Koszul ring.  Our methods will closely mirror those of \cite{2}.


\begin{proposition}\label{main1} We have 
\begin{enumerate}
\item[\bf{(A)}] $M^{(h),L}$ is globally generated for each $h \geq 0$, and
\item[\bf{(B)}] $H^i(M^{(h),L} \otimes B^{\otimes s-i}) = 0$, for all $h \geq 0, ~ s \geq 0,$ and $i \geq 1$.
\end{enumerate}
\end{proposition}
\begin{proof}
We prove both assertions simultaneously by induction on $h$.

First suppose that $h = 0$. $M^{(0),L} = L$ is globally generated because $B$ is. 

Further, for any $s \geq 0$ and $i \geq 1$,
$H^i(M^{(0),L} \otimes B^{\otimes s-i}) = H^i(L \otimes B^{\otimes s-i}) = 
H^i(B^{\otimes r+s-i})= 0$, by (\ref{cm}).

Now fix some $h_1 > 0$ and suppose that the statements {\bf{(A)}} and {\bf{(B)}}
hold for all $h < h_1$. 

So $ M^{(h_1-1),L}$ is globally generated and $M^{(h_1),L}$ is defined. 


We claim that  the multiplication map 
\begin{eqnarray}\label{tmp}
H^0(M^{(h_1),L} \otimes B^{\otimes n}) \otimes H^0(B) \rightarrow H^0( M^{(h_1),L} \otimes B^{\otimes n+1})
\end{eqnarray}
is surjective for all $n \geq 0$. 

By Lemma \ref{key}, this follows if 
\begin{eqnarray}\label{tmp1}
H^i(M^{(h_1),L} \otimes B^{-i}) = 0 ~\text{for all}~ i \geq 1.
\end{eqnarray}

We will first prove (\ref{tmp1}) for $i =1$.

Tensor the sequence (2) corresponding to $F = M^{(h_1-1),L}$ by $B^{\otimes r-1}$. We obtain
$$0 \rightarrow M_{M^{(h_1-1),L}} \otimes B^{\otimes r-1} \rightarrow H^0(M^{(h_1-1),L}) \otimes B^{\otimes r-1} 
\rightarrow M^{(h_1-1),L} \otimes B^{\otimes r-1}\rightarrow 0$$

Taking global sections, we get
\begin{displaymath}
\xymatrix @R=1pc @C=1.5pc{
0 \ar[r]& H^0(M_{M^{(h_1-1),L}} \otimes B^{\otimes r-1}) \ar[r]& H^0(M^{(h_1-1),L}) \otimes H^0(B^{\otimes r-1}) 
\ar[r]^{\gamma} & H^0(M^{(h_1-1),L} \otimes B^{\otimes r-1}) \\
\ar[r] &H^1(M_{M^{(h_1-1),L}} \otimes B^{\otimes r-1})
\ar[r] & H^0(M^{(h_1-1),L}) \otimes H^1(B^{\otimes r-1}) \ar[r] &...
}
\end{displaymath}

By (\ref{cm}), $H^1(B^{\otimes r-1}) = 0$. So $H^1(M_{M^{(h_1-1),L}} \otimes B^{\otimes r-1}) = 0$ 
if and only if $\gamma$ is surjective. 

Now, by Lemmas \ref{horace} and \ref{key}, $\gamma$ is surjective if 
$H^i(M^{(h_1-1),L} \otimes B^{-i}) = 0$ for all $i \geq 1$. 

But this follows from induction hypothesis applied to 
$h_1-1$ and $s=0$. 

Thus $H^1(M_{M^{(h_1-1),L}} \otimes B^{\otimes r-1}) = H^1(M^{(h_1),L} \otimes B^{-1}) = 0$, which 
is the statement (\ref{tmp1}) for $i=1$. 

Now suppose that $i \geq 2$.

Tensor the sequence (2) corresponding to $F = M^{(h_1-1),L}$ by $B^{\otimes r-i}$. 
 We obtain
$$0 \rightarrow M_{M^{(h_1-1),L}} \otimes B^{\otimes r-i} \rightarrow H^0(M^{(h_1-1),L}) \otimes B^{\otimes r-i} 
\rightarrow M^{(h_1-1),L} \otimes B^{\otimes r-i}\rightarrow 0$$

Taking global sections, we get 
\begin{displaymath}
\xymatrix @R=1pc @C=1.5pc{
  ....\ar[r]& 
H^0(M^{(h_1-1),L}) \otimes H^{i-1}(B^{\otimes r-i}) 
\ar[r]& H^{i-1}(M^{(h_1-1),L} \otimes B^{\otimes r-i})\\ 
\ar[r]& H^i(M_{M^{(h_1-1),L}} \otimes B^{\otimes r-i})
\ar[r]& H^0(M^{(h_1-1),L}) \otimes H^i(B^{\otimes r-i}) \ar[r]& ....
}
\end{displaymath}

$H^i(B^{\otimes r-i}) = 0$ by (\ref{cm}).

$H^{i-1}(M^{(h_1-1),L} \otimes B^{\otimes r-i}) = 0$, by induction hypothesis (more precisely, statement {\bf{B}} for $h_1-1$).

Hence $H^i(M_{M^{(h_1-1),L}} \otimes B^{\otimes r-i}) = H^i(M^{(h_1),L} \otimes B^{-i}) = 
0$ for all $i \geq 2$.

This proves (\ref{tmp1}) and hence (\ref{tmp}) for all $i \geq 1$. So we have {\bf{(A)}} by Lemma \ref{globalgen}.

Next, we prove {\bf{(B)}} for $h_1$, $i =1 $ and any $s \geq 0$. 

Tensor the sequence (2) corresponding to $F = M^{(h_1-1),L}$ by $B^{\otimes r+s-1}$. We obtain
$$0 \rightarrow M_{M^{(h_1-1),L}} \otimes B^{\otimes r+s-1} \rightarrow H^0(M^{(h_1-1),L}) \otimes B^{\otimes r+s-1} 
\rightarrow M^{(h_1-1),L} \otimes B^{\otimes r+s-1}\rightarrow 0$$

Taking global sections, we get 
\begin{displaymath}
\xymatrix @R=1pc @C=1.5pc{
0 \ar[r]& H^0(M_{M^{(h_1-1),L}} \otimes B^{\otimes r+s-1}) \ar[r]&H^0(M^{(h_1-1),L}) \otimes H^0(B^{\otimes r+s-1}) 
\ar[r]^{\gamma}& H^0(M^{(h_1-1),L} \otimes B^{\otimes r+s-1})\\ 
\ar[r]& H^1(M_{M^{(h_1-1),L}} \otimes B^{\otimes r+s-1})
\ar[r]& H^0(M^{(h_1-1),L}) \otimes H^1(B^{\otimes r+s-1}) \ar[r]& ....
}
\end{displaymath}

By (\ref{cm}), $H^1(B^{\otimes r+s-1}) = 0$. So $H^1(M_{M^{(h_1-1),L}} \otimes B^{\otimes r+s-1}) = 0$ 
if and only if $\gamma$ is surjective. 

Now, again by Lemmas \ref{horace} and \ref{key}, $\gamma$ is surjective if 
$H^i(M^{(h_1-1),L} \otimes B^{-i}) = 0$ for all $i \geq 1$. 

But this follows from induction hypothesis applied to 
$h_1-1$ and $s=0$. 

Thus $H^1(M_{M^{(h_1-1),L}} \otimes B^{\otimes r+s-1}) = H^1(M^{(h_1),L} \otimes B^{\otimes s-1})= 0$. 

This proves {\bf{(A)}} for $h_1$ and $i=1$.  

Now suppose that $i \geq 2$ and $s \geq 0$.

Tensor the sequence (2) corresponding to $F = M^{(h-1),L}$ by $B^{\otimes r+s-i}$. We obtain
$$0 \rightarrow M_{M^{(h-1),L}} \otimes B^{\otimes r+s-i} \rightarrow H^0(M^{(h-1),L}) \otimes B^{\otimes r+s-i} 
\rightarrow M^{(h-1),L} \otimes B^{\otimes r+s-i}\rightarrow 0$$

Taking global sections, we get
\begin{displaymath}
\xymatrix @R=1pc @C=1.5pc{
... \ar[r]& 
H^0(M^{(h_1-1),L}) \otimes H^{i-1}(B^{\otimes r+s-i}) 
\ar[r]& H^{i-1}(M^{(h_1-1),L} \otimes B^{\otimes r+s-i}) \\
\ar[r]& H^i(M_{M^{(h_1-1),L}} \otimes B^{\otimes r+s-i})
\ar[r]& H^0(M^{(h_1-1),L}) \otimes H^i(B^{\otimes r+s-i}) \ar[]...
}\end{displaymath}
$ H^i(B^{\otimes r+s-i}) = 0$, by (\ref{cm}).

$H^{i-1}(M^{(h_1-1),L} \otimes B^{\otimes r+s-i}) = H^{i-1}(M^{(h_1-1),L} \otimes B^{\otimes r+s-1-(i-1)}) = 
0$, by induction hypothesis because $r+s-1 \geq s \geq 0$.

So $H^i(M_{M^{(h_1-1),L}} \otimes B^{\otimes r+s-i}) = H^i(M^{(h_1),L}  \otimes B^{\otimes s-i}) = 0$, as required.
\end{proof}

\begin{theorem}\label{final}
Let $X$ be a projective variety over a field $k$. Let $B$ a base point free ample bundle on $X$ with reg($B) =r-1$. Let $L = B^{\otimes n}$ with $n \geq r$.  Then $R(L)$ is a Koszul ring.
\end{theorem}
\begin{proof}
Since $B$ is $(n-1)$-regular, Proposition \ref{main1} implies that, $H^1(M^{(h),L} \otimes B^{\otimes s-1}) = 0$ for all $s \geq 0$ and $h \geq 0$.

So for $s \geq 0$ and $h \geq 0$,
$H^1(M^{(h),L} \otimes L^s) = H^1(M^{(h),L} \otimes B^{\otimes rs}) = 0$. 

By Proposition \ref{p}, it follows that $R(L)$ is a Koszul ring. 
\end{proof}
\section{Minimal irregular surfaces of general type}\label{part2}
Let $X$ be a nonsingular projective minimal\footnote{A surface $X$ is {\it minimal} if every birational morphism $X \rightarrow Y$ is an isomorphism.} irregular\footnote{The {\it irregularity} $q$ of a surface $X$ over a field $k$ is defined to be  $q = dim_k ~H^1(\ox)$. We say that $X$ is irregular if $q > 0$.} surface of general type\footnote{A surface $X$ is of {\it general type} if its Kodaira dimension $\kappa(X) = 2$.} over the complex number field
$\mathbb{C}$. Let $K_X$ be the canonical line bundle on $X$. Suppose that $K_X$ is base point free. 

{\bf Notation:} 
We write $L \equiv L'$ if the line bundles $L$ and $L'$ are numerically equivalent. We write 
$L \cdot L'$ to denote the intersection number of $L$ and $L'$.

Let $B$ be a base point free and ample divisor on $X$ such that $B'$ is base point free for all $B'\equiv B$ and $H^1(B') = 0$. Assume that 
$B^2 > B \cdot K_X$.


Let $L = K_X \otimes B^{\otimes n}$, where $n \geq 2$. 

Set 
$\begin{displaystyle}R(L) = \bigoplus_{n=0}^{\infty} H^0(L^{\otimes n})
\end{displaystyle}.$ Our goal is to prove that $R(L)$ is Koszul. Our proof is similar to proofs in \cite{4}. Theorem 5.14 in \cite{2} proves an analogous result for regular surfaces.

\subsection{Required lemmas}

In this subsection we will prove some lemmas that will be used later in the proof of the main theorem. 

A divisor $D$ on $X$ is {\it nef} if $D \cdot C \geq 0$ for every irreducible curve $C$ in $X$. 
$D$ is said to be {\it big} if a multiple $mD$, $m \in \mathbb{N}.$  defines a birational map of $X$ to a projective space.

As $X$ is a minimal surface of general type, $K_X$ is nef and big. In fact, a surface is minimal of general type if and only if $K_X$ is nef and big.

\begin{lemma}\label{kv}
{\bf \text (Kawamata - Viehweg vanishing)} Let $X$ be a nonsingular projective variety over the complex number field $\mathbb{C}$. Let $D$ be a nef and big divisor on $X$. Then $$H^i(K_X \otimes D) = 0, \text {~for~} i > 0.$$ 
\end{lemma}

For a proof, see \cite{6} or \cite{7}. We will refer to this result simply as K-V vanishing.

Recall that ${Pic}^0(X)$ denotes the group of divisors on $X$ which are algebraically equivalent to zero modulo linear equivalence. 
\begin{lemma}\label{pic}
There exists a divisor $E \in {Pic}^0(X)$ such that $E^{\otimes 2} \neq \ox$.
\end{lemma}
\begin{proof} We have the exponential sequence 
$$0 \rightarrow \mathbb{Z} \rightarrow \oxh \rightarrow {\cal O}^{\star}_{X_h} \rightarrow 0,$$
where $X_h$ is the complex analytic space associated to $X$. 

Consider the resulting long exact sequence in cohomologies. Applying Serre's GAGA and identifying the $Pic(X)$ with $H^1(X,{\cal O}^{\star}_{X})$, we obtain an exact sequence 
$$0 \rightarrow H^1(X_h, \mathbb{Z}) \rightarrow H^1(X, \ox) \rightarrow Pic(X) \rightarrow H^2(X_h, \mathbb{Z}) \rightarrow H^2(X,\ox) \rightarrow ...$$

This gives ${Pic}^0(X) \cong H^1(X,\ox)/H^1(X_h,\mathbb{Z})$. This is an abelian variety.

For more details on this see the discussion in Appendix B.5 in \cite{9}.

Since $X$ is irregular $H^1(\ox) \neq 0$. So ${Pic}^0(X)$ is a nontrivial abelian variety and hence contains 2-torsion elements. 
\end{proof}
\begin{lemma}\label{-10} 
$H^1(B_1 \otimes B_2 \otimes \otimes ... \otimes B_n) = 0$, for line bundles 
$B_1 \equiv B_2 \equiv ... \equiv B_n \equiv B$ and $n \geq 1$.
\end{lemma}
\begin{proof}
Let $C \in |B|$ be a smooth curve. We have for every $i$, 
$$\text {deg}(B_i \otimes \oc) = B_i \cdot C = B^2$$ 


If $n > 3$, deg$(B_1 \otimes ... \otimes B_n \otimes \oc) = nB^2 > 2B^2 = B^2 + B^2 \geq B^2 +B \cdot K_X = 2g(C) -2$, where $g(C)$ is the genus of $C$. So 
$H^1(B_1 \otimes ... \otimes B_n \otimes \oc) = 0$.

We have the short exact sequence: 
\begin{eqnarray}\label{ses}
0 \rightarrow B^{-1}  \rightarrow {\cal O}_X \rightarrow {\cal O}_C \rightarrow 0 
\end{eqnarray}
Tensoring with $B_1  \otimes  B_2 \otimes  B_3$, we get
$$0 \rightarrow B_1 \otimes B_2 \otimes B_3 \otimes B^{-1} \rightarrow B_1  \otimes  B_2 \otimes  B_3 \rightarrow B_1  \otimes  B_2 \otimes  B_3 \otimes {\cal O}_C \rightarrow 0.$$

Note that $B_3 \otimes B^{-1}  \equiv {\cal O}_B$, so we can write $B_1  \otimes  B_2 \otimes  B_3=B_1 \otimes {B_2}'$, where $B_2' \equiv B_2 \equiv B$.

Taking the long exact sequence in cohomology of the above short exact sequence we obtain,
$$H^1(B_1 \otimes B_2') \rightarrow H^1(B_1 \otimes B_2 \otimes B_3) \rightarrow 
H^1(B_1 \otimes B_2 \otimes B_3 \otimes \oc).$$

Since $H^1(B_1 \otimes B_2\otimes B_3 \otimes \oc) = 0$, it is enough to prove
the theorem for $n = 2$.

Exactly as above, we have the following exact sequence
$$H^1(B_1) \rightarrow H^1(B_1 \otimes B_2) \rightarrow 
H^1(B_1 \otimes B_2 \otimes \oc).$$

$H^1(B_1)=0$ by hypothesis. It is enough to prove that $H^1(B_1 \otimes B_2 \otimes \oc) = 0$.

As before, deg$(B_1 \otimes B_2 \otimes \oc) = 2B^2 \geq B^2+ B.K_X = 2g(C)-2$. 

If deg$(B_1 \otimes B_2 \otimes \oc) > 2g(C)-2$, then we are done. 

Suppose that  
deg$(B_1 \otimes B_2 \otimes \oc) = 2g(C)-2$. Note that this implies that $B^2 = B.K_X$.

If $B_1 \otimes B_2 \otimes \oc \neq K_C$, then $H^1(B_1 \otimes B_2 \otimes \oc) = 0$ and we are done. 
Here $K_C$ denotes the canonical divisor of $C$.

Assume that $B_1 \otimes B_2 \otimes \oc = K_C$. By adjunction, we have $K_C =  K_X \otimes B \otimes \oc$. So 
$B_1 \otimes B_2 \otimes \oc = K_X \otimes B \otimes \oc$. This gives
$B_1 \otimes B_2 \otimes K^{-1}_X \otimes B^{-1} \otimes \oc = {\cal O}_C$.

Tensoring (\ref{ses}) with $B_1 \otimes B_2 \otimes {K}^{-1}_X \otimes B^{-1}$, we obtain
$$0 \rightarrow B_1 \otimes B_2 \otimes {K}^{-1}_X \otimes B^{-2}  \rightarrow B_1 \otimes B_2 \otimes {K}^{-1}_X \otimes B^{-1}  \rightarrow 
B_1 \otimes B_2 \otimes {K}^{-1}_X \otimes B^{-1} \otimes {\cal O}_C \rightarrow 0.$$

Taking cohomology long exact sequence, we have
\begin{displaymath}
\xymatrix @R=1pc @C=1.5pc{
H^0(B_1 \otimes B_2 \otimes K^{-1}_X \otimes B^{\otimes -2}) \ar[r]& 
H^0(B_1 \otimes B_2 \otimes K^{-1}_X \otimes B^{-1}) \ar[r]&\\
H^0(B_1 \otimes B_2 \otimes K^{-1}_X \otimes B^{-1} \otimes \oc) 
\ar[r]&
H^1(B_1 \otimes B_2 \otimes K^{-1}_X \otimes B^{\otimes -2})
}\end{displaymath}
Now $H^0(B_1 \otimes B_2 \otimes K^{-1}_X \otimes B^{\otimes -2})  = 
H^2(K^{\otimes 2}_X \otimes {B_1}^{-1} \otimes {B_2}^{-1} \otimes B^{\otimes 2})$, by Serre duality. Since 
${B_1}^{-1} \otimes {B_2}^{-1} \otimes B^{\otimes 2} \equiv \ox$ and $K_X$ is nef and big ($X$ is minimal of general type), it follows that $H^2(K^{\otimes 2}_X \otimes {B_1}^{-1} \otimes {B_2}^{-1} \otimes B^{\otimes 2}) = 0$ by K-V vanishing.

Similarly, $H^1(B_1 \otimes B_2 \otimes K^{-1}_X \otimes B^{\otimes -2})  = 
H^1(K^{\otimes 2}_X \otimes {B_1}^{-1} \otimes {B_2}^{-1} \otimes B^{\otimes 2}) = 0$. 


Thus we obtain 

$H^0(B_1 \otimes B_2 \otimes K^{-1}_X \otimes B^{-1}) \cong 
H^0(B_1 \otimes B_2 \otimes K^{-1}_X \otimes B^{-1} \otimes \oc) \cong H^0(\oc) \cong k$.

So $H^0(B_1 \otimes B_2 \otimes K^{-1}_X \otimes B^{-1}) \neq 0$ and 
$B_1 \otimes B_2 \otimes K^{-1}_X \otimes B^{-1}$ is effective. 
But $$B\cdot (B_1 \otimes B_2 \otimes K^{-1}_X \otimes B^{-1}) = B \cdot B_1 + B \cdot B_2 -B \cdot K_X - B^2 = B^2 - B \cdot K_X = 0$$

So $B_1 \otimes B_2 \otimes K^{-1}_X \otimes B^{-1} = \ox \Rightarrow B_1 \otimes B_2 = 
K_X \otimes B$. Finally, $H^1(B_1 \otimes B_2) = H^1(K_X \otimes B) = 0$, by K-V vanishing, thus concluding the proof. \end{proof}

\begin{lemma}\label{0}
$H^2(B^{\otimes n} \otimes \delta) = 0$ for $n \geq 1$ and any numerically trivial line bundle $\delta$.
\end{lemma}
\begin{proof}
$H^2(B^{\otimes n} \otimes \delta) = H^0(K_X \otimes B^{-n} \otimes {\delta}^{-1})$, by Serre duality.

If $H^0(K_X \otimes B^{-n} \otimes {\delta}^{-1}) \neq 0$, then there is an effective divisor $D$ that is linearly equivalent to $K_X- nB -\delta$. (By abuse of notation, we denote the divisor associated to a line bundle by the same letter.) So we have $B \cdot D = B \cdot (K_X-nB-\delta) \geq 0$, because $B$ is ample. So $B \cdot K_X \geq  n B^2$. But this contradicts the hypothesis that 
$B^2 >  B \cdot K_X$.
\end{proof}
\begin{lemma}\label{1}
$H^2(K^{-1}_X \otimes B^{\otimes n} \otimes \delta) = 0$, for $n \geq 2$. 
\end{lemma}
\begin{proof}
By Serre duality, $H^2(K^{-1}_X \otimes B^{\otimes n} \otimes \delta) = H^0({K}^{\otimes 2}_X \otimes B^{-n} \otimes {\delta}^{-1})$.

If $H^0({K}^{\otimes 2}_X \otimes B^{-n} \otimes {\delta}^{-1}) \neq 0$, then there is an effective divisor $D$ that is linearly equivalent to $2K_X - nB - \delta$. So we have $B \cdot D = B \cdot  (2K_X-nB - \delta) \geq 0$, because $B$ is ample. So $2B \cdot K_X \geq n B^2$. But this contradicts the hypothesis that $B^2 > B \cdot K_X$.
\end{proof} 

\subsection{Main theorem}
In this subsection we will prove our main theorem: $R(L)$ is a Koszul ring. 


By Lemma \ref{pic}, there exists $E \in Pic^0(X)$ such that $E^{\otimes 2} \neq \ox$. Note that $E$ is numerically trivial.

Set $B_1 = B \otimes E^{-1}$ and $B_2 = B \otimes E$.  Then $L = B_1 \otimes B_2 \otimes K_X$. 

Let $\delta$ be a numerically trivial line bundle such that $${\delta}^{\otimes 2} \neq E^{\otimes 2}.$$ Let $m, r, t$ be non-negative integers such that \underline{$m+r+t >0$}.

Recall the definition of $H^1(M^{(h),L})$ from Section \ref{section1}. $H^1(M^{(h-1),L})$ has to be globally generated to define  $H^1(M^{(h),L})$. The discussion below will establish that $H^1(M^{(h),L})$  is globally generated for all $h \geq 0$. 

Consider the following statements for a non-negative integer $h$.\\

{\rm \bf (V$_h$)}
~~~~~~~~~~~~$H^1(M^{(h),L} \otimes {B_1}^{\otimes m} {B_2}^{\otimes r} \otimes K^{\otimes t} \otimes \delta) = 0 $\\\\

{\rm \bf (S$_h$)}
$  \left\{
\begin{array}{l}
{\rm The ~multiplication ~map~}\\
H^0(M^{(h),L}) \otimes H^0({B_1}^{\otimes m} \otimes {B_2}^{\otimes r} \otimes K^{\otimes t} \otimes \delta) \rightarrow H^0(M^{(h),L} \otimes {B_1}^{\otimes m} \otimes {B_2}^{\otimes r} \otimes K^{\otimes t} \otimes \delta ) \\ {\rm ~is ~surjective}
\end{array} \right.
$\\\\


Our goal is to prove that {\bf (S$_h$)} and {\bf (V$_h$)} hold for all $h \geq 0$. 

\begin{lemma}\label{h=0}
The statements {\bf (S$_0$)} and {\bf (V$_0$)} hold.
\end{lemma}
\begin{proof}
{\bf (V$_0$)}: $H^1(L \otimes {B_1}^{\otimes m} {B_2}^{\otimes r} \otimes K^{\otimes t} \otimes \delta) = 0$ holds by K-V vanishing.

To prove {\bf (S$_0$)}, we will use Lemma \ref{horace} iteratively.  First, let us observe that the following map is surjective for $m \geq 0$:
\begin{eqnarray}\label{0a}
H^0(L \otimes {B_1}^{\otimes m}) \otimes H^0({B_1} \otimes \delta) \rightarrow H^0(L \otimes {B_1}^{\otimes m+1} \otimes \delta).
\end{eqnarray}

By Lemma \ref{key}, we need $H^1(L \otimes {B_1}^{\otimes m-1} \otimes {\delta}^{-1}) = 0$ and $H^2(L \otimes {B_1}^{\otimes m-2} \otimes {\delta}^{-2}) = 0$.

$H^1(L \otimes {B_1}^{\otimes m-1} \otimes {\delta}^{-1}) = H^1(K_X  \otimes {B_1}^{\otimes m} \otimes B_2 \otimes {\delta}^{-1}) = 0$, by K-V vanishing. 

$H^2(L \otimes {B_1}^{\otimes m-2} \otimes {\delta}^{-2}) = H^2(K_X  \otimes {B_1}^{\otimes m-1} \otimes B_2 \otimes {\delta}^{-2})$. 

If $m=0$, then $H^2(K_X  \otimes {B_1}^{-1} \otimes B_2 \otimes {\delta}^{-2}) = H^2(K_X \otimes E^{\otimes 2} \otimes {\delta}^{ -2}) = H^0(E^{-2} \otimes {\delta}^{\otimes 2}) = 0$ because 
$E^{-2} \otimes {\delta}^{\otimes 2} \neq \ox$ is numerically trivial. 

If $m > 0$, then $H^2(K_X  \otimes {B_1}^{\otimes m-1} \otimes B_2 \otimes {\delta}^{-2}) = H^2(K_X \otimes B^{\otimes m} \otimes \delta_1) = 0$ by K-V vanishing ($\delta_1$ is a numerically trivial line bundle).

Second, let us show that the following map is surjective for $m, r \geq 0$:
\begin{eqnarray}\label{0b}
H^0(L \otimes {B_1}^{\otimes m} \otimes {B_2}^{\otimes r}) \otimes H^0(B_2  \otimes \delta) \rightarrow H^0(L \otimes {B_1}^{\otimes m} \otimes {B_2}^{\otimes r+1}  \otimes \delta).
\end{eqnarray}
By Lemma \ref{key}, we need $H^1(L \otimes {B_1}^{\otimes m} \otimes {B_2}^{\otimes r-1}  \otimes {\delta}^{-1}) = 0$ and $H^2(L \otimes {B_1}^{\otimes m}  \otimes {B_2}^{\otimes r-2}  \otimes {\delta}^{-2}) = 0$.

$H^1(L \otimes {B_1}^{\otimes m}  \otimes {B_2}^{\otimes r-1}  \otimes {\delta}^{-1}) = H^1(K_X  \otimes {B_1}^{\otimes m+1} \otimes {B_2}^{\otimes r}  \otimes {\delta}^{-1}) = 0$, by K-V vanishing. 

$H^2(L \otimes {B_1}^{\otimes m}  \otimes {B_2}^{\otimes r-2}  \otimes {\delta}^{-2}) = H^2(K_X  \otimes {B_1}^{\otimes m+1} \otimes {B_2}^{\otimes r-1}  \otimes {\delta}^{-2})$. 

If $m=r=0$, then $H^2(K_X  \otimes {B_1}^{\otimes m+1} \otimes {B_2}^{-1}  \otimes {\delta}^{-2}) = H^2(K_X \otimes E^{\otimes 2} \otimes {\delta}^{ -2}) = H^0(E^{-2} \otimes {\delta}^{\otimes 2}) = 0$ because as above 
$E^{-2} \otimes {\delta}^{\otimes 2} \neq \ox$ is numerically trivial. 

If $m+r > 0$, then $H^2(K_X  \otimes {B_1}^{\otimes m+1} \otimes {B_2}^{\otimes r-1}  \otimes {\delta}^{-2}) = H^2(K_X \otimes B^{\otimes m+r} \otimes \delta_1) = 0$ by K-V vanishing ($\delta_1$ is a numerically trivial line bundle).

Finally, we will prove the following map is surjective:
\begin{eqnarray}\label{0c}
H^0(L \otimes {B_1}^{\otimes m} \otimes {B_2}^{\otimes r} \otimes {K_X}^{\otimes t}) \otimes H^0(K_X  \otimes \delta) \rightarrow H^0(L \otimes {B_1}^{\otimes m} \otimes {B_2}^{\otimes r}  \otimes {K_X}^{\otimes t+1}  \otimes \delta).
\end{eqnarray}

By Lemma \ref{key}, we need $H^1(L \otimes {B_1}^{\otimes m} \otimes {B_2}^{\otimes r}  \otimes {K_X}^{\otimes t-1} \otimes {\delta}^{-1}) = 0$ and $H^2(L \otimes {B_1}^{\otimes m}  \otimes {B_2}^{\otimes r} \otimes {K_X}^{\otimes t-2} \otimes {\delta}^{-2}) = 0$.

$H^1(L \otimes {B_1}^{\otimes m}  \otimes {B_2}^{\otimes r}  \otimes {K_X}^{\otimes t-1} \otimes {\delta}^{-1}) = H^1( {K_X}^{\otimes t}  \otimes {B_1}^{\otimes m+1} \otimes {B_2}^{\otimes r+1}  \otimes {\delta}^{-1}) $. If $t=0$, this is zero by Lemma \ref{-10}. If $t>0$, it is zero by K-V vanishing. 

$H^2(L \otimes {B_1}^{\otimes m}  \otimes {B_2}^{\otimes r}  \otimes {K_X}^{\otimes t-2} \otimes {\delta}^{-2}) = H^2({K_X}^{\otimes t-1}  \otimes {B_1}^{\otimes m+1} \otimes {B_2}^{\otimes r+1}  \otimes {\delta}^{-2})$. 

If $t=0$, then this $H^2$ is zero by Lemma \ref{1}. If $t=1$ it is zero by Lemma \ref{-10}. If $t>0$, it is zero by K-V vanishing. 

The proof is now complete by the surjectivity  of (\ref{0a}), (\ref{0b}) and (\ref{0c}), and Lemma \ref{horace}.
\end{proof}

Note that (\ref{0a}) shows that $L$ is globally generated (taking $\delta = \ox$ and applying Lemma \ref{globalgen}). So we can define $M^{(1),L}$.

\begin{theorem}\label{main}
The statements {\bf (V$_h$)} and {\bf (S$_h$)} hold all $h \geq 0$.

\end{theorem}

\begin{proof} The proof is by induction on $h$. Both statements hold when $h = 0$ by Lemma \ref{h=0}.

Suppose that the statements hold for all nonnegative integers $\leq h-1$ for some $h \geq 1$.

Proving that {\bf (V$_h$)} holds is easy:
by Remark \ref{-1}, {\bf (V$_h$)} follows if
\begin{itemize}
\item [(i)] $H^0(M^{(h-1),L}) \otimes H^0(L \otimes {B_1}^{\otimes m} \otimes {B_2}^{\otimes r} \otimes K^{\otimes t} \otimes \delta) \rightarrow H^0(M^{(h-1),L} \otimes L \otimes {B_1}^{\otimes m} \otimes {B_2}^{\otimes r} \otimes K^{\otimes t} \otimes \delta )$
is surjective, and
\item [(ii)] $H^1({B_1}^{\otimes m} \otimes {B_2}^{\otimes r} \otimes K^{\otimes t} \otimes \delta) = 0$.
\end{itemize}

(ii) follows by K-V vanishing. (i) is simply the statement 
{\bf (S$_{h-1}$)}.

To prove {\bf (S$_{h}$)}, we will need to do some work. We are going to use Lemma \ref{horace} iteratively. Lemma \ref{horace} allows us to prove the surjectivity separately for $B_1$, $B_2$ and 
$K_X$, as in Lemma \ref{h=0}.  We will deal with these three cases in the three lemmas that follow.

First, we prove the following
\begin{lemma}\label{ma}
The multiplication map 
$$H^0(M^{(h),L} \otimes {B_1}^{\otimes m}) \otimes H^0({B_1} \otimes \delta) \rightarrow H^0(M^{(h),L} \otimes {B_1}^{\otimes m+1} \otimes \delta)$$ is surjective for $m \geq 0$.
\end{lemma}
\begin{proof}
We use Lemma \ref{key}. We need the following two statements for $m \geq 0$:
\begin{eqnarray}
H^1(M^{(h),L} \otimes {B_1}^{\otimes m-1} \otimes {\delta}^{-1}) = 
H^1(M_{M^{(h-1),L}} \otimes L \otimes {B_1}^{\otimes m-1} \otimes {\delta}^{-1}) = 0\label{ma1}\\
H^2((M^{(h),L} \otimes {B_1}^{\otimes m-2} \otimes {\delta}^{-2}) = 
H^2(M_{M^{(h-1),L}} \otimes L \otimes {B_1}^{\otimes m-2} \otimes {\delta}^{-2}) = 0\label{ma2}
\end{eqnarray}\\

By Remark \ref{-1}, (\ref{ma1}) follows if 
$H^1( L \otimes {B_1}^{\otimes m-1}) = 0$ and if 
the following map is surjective:
$$H^0(M^{(h-1),L}) \otimes H^0(L \otimes {B_1}^{\otimes m-1}) 
\rightarrow H^0(M^{(h-1),L} \otimes L \otimes {B_1}^{\otimes m-1}).$$

The $H^1$ is zero by K-V vanishing and the surjectivity is simply {\bf (S$_{h-1})$}.\\

By Remark \ref{-2}, 
(\ref{ma2}) follows if 
$H^1(M^{(h-1),L} \otimes L \otimes {B_1}^{\otimes m-2} \otimes {\delta}^{-2}) = 0$ 
and 

$H^2(L \otimes {B_1}^{\otimes m-2} \otimes {\delta}^{-2} ) = 0$

The $H^1$ vanishes by {\bf (V$_{h-1}$)}. If $m=0$, 
$H^2(K_X \otimes E^{\otimes 2} \otimes {\delta}^{ -2}) = 0$, as in the proof of Lemma \ref{h=0}.
$H^2$ is zero by K-V vanishing. If $m > 0$, then $H^2$ is zero by K-V vanishing. 
\end{proof}


Lemma \ref{ma} implies that $M^{(h),L}$ is globally generated (by Lemma \ref{globalgen}) for all $h \geq 0$.

Now to the next step:

\begin{lemma}\label{mb}
The multiplication map 
$$H^0(M^{(h),L} \otimes {B_1}^{\otimes m} \otimes {B_2}^{\otimes r}) \otimes H^0({B_2} \otimes \delta) \rightarrow H^0(M^{(h),L} \otimes {B_1}^{\otimes m} \otimes {B_2}^{\otimes r+1} \otimes \delta)$$ is surjective for $m, r \geq 0$.
\end{lemma}
\begin{proof}
According to Lemma \ref{key}, we need the following two statements for $m, r \geq 0$:
\begin{eqnarray}
H^1(M_{M^{(h-1),L}} \otimes L \otimes {B_1}^{\otimes m} \otimes {B_2}^{\otimes r-1} \otimes {\delta}^{-1}) = 0\label{mb1}\\
H^2(M_{M^{(h-1),L}} \otimes L \otimes {B_1}^{\otimes m} \otimes {B_2}^{\otimes r-2} \otimes {\delta}^{-2}) = 0\label{mb2}
\end{eqnarray}

By Remark \ref{-1}, 
(\ref{mb1}) follows if $H^1(L \otimes {B_1}^{\otimes m} \otimes {B_2}^{r-1} \otimes {\delta}^{-1}) = 0$ and if 
the following map is surjective:
$$H^0(M^{(h-1),L}) \otimes H^0(L \otimes {B_1}^{\otimes m} \otimes {B_2}^{r-1} \otimes {\delta}^{-1}) 
\rightarrow H^0(M^{(h-1),L} \otimes L \otimes {B_1}^{\otimes m} \otimes {B_2}^{r-1} \otimes {\delta}^{-1}).$$

The $H^1$ is zero by K-V vanishing and the surjectivity is simply {\bf (S$_{h-1})$}.\\

By Remark \ref{-2}, 
(\ref{mb2}) follows if 
$H^1(M^{(h-1),L} \otimes L \otimes {B_1}^{\otimes m} \otimes {B_2}^{\otimes r-2} \otimes {\delta}^{-2}) = 0$ 
and 

$H^2(L \otimes {B_1}^{\otimes m} \otimes {B_2}^{\otimes r-2} \otimes {\delta}^{-2} ) = 0$

The $H^1$ vanishes by {\bf (V$_{h-1}$)}. If $m=r=0$, then $H^2$ is zero as in Lemma \ref{mb}. If 
$m+r > 0$, then $H^2$ is zero by K-V vanishing.
\end{proof}

Finally we have the following 
\begin{lemma}\label{mc}
The multiplication map 
$$H^0(M^{(h),L} \otimes {B_1}^{\otimes m} \otimes {B_2}^{\otimes r} \otimes K^{\otimes t}) \otimes H^0({K \otimes \delta}) \rightarrow H^0(M^{(h),L} \otimes {B_1}^{\otimes m} \otimes {B_2}^{\otimes r}  \otimes K^{\otimes t+1} \otimes \delta)$$ is surjective for $m, r, t \geq 0$.
\end{lemma}
\begin{proof}
According to Lemma \ref{key}, we need the following two statements for $m, r, t \geq 0$:
\begin{eqnarray}
H^1(M_{M^{(h-1),L}} \otimes L \otimes {B_1}^{\otimes m} \otimes {B_2}^{\otimes r}  \otimes K^{\otimes t-1} \otimes  {\delta}^{-1}) = 0\label{mc1}\\
H^2(M_{M^{(h-1),L}} \otimes L \otimes {B_1}^{\otimes m} \otimes {B_2}^{\otimes r}  \otimes K^{\otimes t-2} \otimes {\delta}^{-2}) = 0\label{mc2}
\end{eqnarray}

By Remark \ref{-1}, 
(\ref{mc1}) follows if $H^1(L \otimes {B_1}^{\otimes m} \otimes {B_2}^{\otimes r}  \otimes K^{\otimes t-1}\otimes {\delta}^{-1}) = 0$ and if 
the following map is surjective:
$$H^0(M^{(h-1),L}) \otimes H^0(L \otimes {B_1}^{\otimes m} \otimes {B_2}^{\otimes r}  \otimes K^{\otimes t-1} \otimes {\delta}^{-1}) 
\rightarrow H^0(M^{(h-1),L} \otimes L \otimes {B_1}^{\otimes m} \otimes {B_2}^{r-1} \otimes {\delta}^{-1}).$$

The surjectivity is simply {\bf (S$_{h-1})$}. For $t \geq 1$, the $H^1$ is zero by K-V vanishing. 
For $t=0$, the vanishing follows from Lemma \ref{0}. \\

By Remark \ref{-2}, 
(\ref{mc2}) follows if 
$H^1(M^{(h-1),L} \otimes L \otimes {B_1}^{\otimes m} \otimes {B_2}^{\otimes r}  \otimes K^{\otimes t-1} \otimes {\delta}^{-2}) = 0$ 
and 

$H^2(L \otimes  {B_1}^{\otimes m} \otimes {B_2}^{\otimes r}  \otimes K^{\otimes t-1}\otimes {\delta}^{-2} ) = 0$

The $H^1$ vanishes by {\bf (V$_{h-1}$)}. For $t \geq 2$, $H^2$ is zero by K-V vanishing; for
$t =1$, $H^2$ vanishes by Lemma \ref{0}; for $t=0$, it vanishes by Lemma \ref{1}.
\end{proof}

{\bf (S$_{h}$)} now from Lemma \ref{horace} and Lemmas \ref{ma}, \ref{mb}, \ref{mc}.
\end{proof}

\begin{theorem}\label{last}
Let $X$ be an nonsingular projective minimal irregular surface of general type over $\mathbb{C}$. Suppose that the canonical divisor $K_X$ of $X$ is base point free. Let $B$ be a base point free ample divisor on $X$ such that $B'$ is base point free for all $B'\equiv B$ and $H^1(B') = 0$. Assume that $B^2 > B \cdot K_X$. Let $L = K_X \otimes B^{\otimes n}$, where $n \geq 2$. Then $R(L)$ is a Koszul ring. 
\end{theorem} 
\begin{proof}
By Theorem \ref{main}, {\bf (V$_h$)} and {\bf (S$_h$)} hold for all $h \geq 0$. 

If $s > 0$, then taking $s=m=r=t$ and $\delta = \ox$, {\bf (V$_h$)} gives us $H^1(M^{(h),L} \otimes L^{\otimes s}) = 0$. 

If $s=0$, we need to prove that $H^1(M^{(h),L}) = 0$ for $h \geq 0$. If $h=0$ this follows by K-V vanishing. Suppose that $h > 0$. We need to prove that $H^1(M^{(h),L}) = H^1(M_{M^{(h-1),L}} \otimes L) = 0$. By Remark \ref{-1}, this follows if the multiplication map
$$H^0(M^{(h-1),L}) \otimes H^0(L) \rightarrow H^0(M^{(h-1),L} \otimes L)$$
is surjective and if $H^1(L) = 0$. Surjectivity is simply the statement {\bf (S$_{h-1}$)}
and $H^1(L) = 0$ by K-V vanishing.

This implies $H^1(M^{(h),L} \otimes L^{\otimes s}) = 0$ for $h \geq 0$ and $s \geq 0$, thereby proving that $R(L)$ is Koszul, by Proposition \ref{p}.
\end{proof}

\bibliographystyle{plain}

\begin{thebibliography}{99}
\bibitem{12} Avramov, L.L.,  and Eisenbud, D., {\it  Regularity of modules over a Koszul algebra},  J. Algebra 153, 85-90 (1992).
\bibitem{16} Backelin, J.,  {\it On the rates of growth of the homologies of Veronese subrings},  Algebra, algebraic topology and their interactions (Stockholm, 
1983), Vol 1183, Lecture Notes in Mathematics, 79-100, Springer-Verlag, 1986. 
\bibitem{13} Beilinson, A.A., Ginsburg V.A., and Schechtman V.V., {\it Koszul duality}, J. Geom. Phys. 5 (1988), 317-350
\bibitem{14} Beilinson, A.A., Ginzburg V.A., and Soergel, W., {\it Koszul duality patterns in representation theory}, J. Amer. Math. Soc. 9 (1996), 473-527. 
\bibitem{17} Butler, David, {\it Normal generation of vector bundles over a curve},  J. Differential Geom., 39, 1-34, 1994. 
\bibitem{23} Eisenbud, D., Reeves, A., and Totaro, B., {\it Initial ideals, Veronese
subrings, and rates of algebras}, Adv. Math., 109, 168-187, 1994.
\bibitem{1} Gallego, Francisco Javier and Purnaprajna, B. P., {\it Normal presentation on elliptic ruled surfaces},  J. Algebra  186  (1996),  no. 2, 597-625.
\bibitem{2} Gallego, F. J. and Purnaprajna, B. P., {\it Projective normality and syzygies of algebraic surfaces},  J. Reine Angew. Math.  506  (1999), 145-180
\bibitem{9} Hartshorne, Robin, {\it Algebraic Geometry}, Graduate Texts in Mathematics, 52, Springer-Verlag, 1977.
\bibitem{22} Hering Milena,  {\it Multigraded regularity and Koszul property}, Preprint, arXiv:0712.2251v1.
\bibitem{11} Herzog, J. and and Iyengar, S.,  {\it Koszul modules}, J. Pure Appl. Algebra 201 (2005), 154-188. 
\bibitem{6} Kawamata, Yujiro, {\it A generalization of Kodaira-Ramanujam's vanishing theorem}, Math. Ann. 261 (1982), no. 1, 43-46.
\bibitem{8} Lazarsfeld, Robert, {\it Positivity in Algebraic Geometry I}, Springer-Verlag, 2004.  
\bibitem{5} Mumford, David, {\it Varieties defined by quadratic equations, Corso CIME in Questions on Algebraic Varieties}, Rome (1970), 30-100.
\bibitem{3} Pareschi, G., {\it Koszul algebras associated to adjunction bundles}, J. of Algebra 157 (1993) 161-169.
\bibitem{20} Pareschi, G. and Purnaprajna, B.P.,  {\it Canonical ring of a curve is Koszul: a simple proof}, Illinois J. Math., 41 (1997), 266-271. 
\bibitem{15} Polishchuk, A. and Positselski, L., {\it Quadratic algebras}, University Lecture Series, Vol 37, American Mathematical Society, Providence, RI, 2005. 
\bibitem{18} Polishchuk A., {\it On the Koszul property of the homogeneous coordinate ring of a curve},  J. Algebra, 178 (1995), 122-135. 
\bibitem{4} Purnaprajna, B.P., {\it Some results on surfaces of general type},  Canad. J. Math.  57  (2005), 4, 724-749.
\bibitem{10} Priddy, Stewart, {\it Koszul resolutions}, Trans. Amer. Math. Soc. 152 (1970), 39-60. 
\bibitem{21} Sturmfels, B., {\it Four counterexamples in combinatorial algebraic geometry}, J. Algebra, 230 (2000), 282-294. 
\bibitem{7} Viehweg, Eckart, {\it Vanishing theorems}, J. Reine Angew. Math. 335 (1982), 1-8.
\bibitem{19} Vishik, A. and Finkelberg, M., {\it The coordinate ring of general curve of genus g $\geq$  5 is Koszul} , J. Algebra 162 (1993), 535-539. 
\end{thebibliography}

Department of Mathematics, University of Kansas, Lawrence, KS 66049\\
Email address: khanuma@math.ku.edu
\end{document}